\documentstyle[12pt]{article}
\normalbaselines
\catcode `\@=11
\def\ot{\otimes}
\def\tg{\triangle}
\def\ra{\longrightarrow}
\def\ca{{\cal A}}

\def\cc{{\cal C}}

\def\cm{{\cal M}}
\def\ps#1,#2,{\Psi_{{#1}{,}{#2}}}
\def\ep{\epsilon}

\def\a{\alpha}
\def\b{\beta}

\def\g{\gamma}

\def\ep{\epsilon}
\def\si{\sigma}

\def\th{\Theta}
\def\ba{\begin{array}}
\def\ea{\end{array}}  
\def\be{\begin{equation}}
\def\ee{\end{equation}}
\newfont{\numb}{msbm10}
\def\com{\mbox{\numb C}}
\def\real{\mbox{\numb R}}
\def\tr{\triangleright}
\def\lb{\langle}
\def\rb{\rangle}
\@addtoreset{equation}{section}

\def\section{\@startsection {section}{1}{\z@}{-3.5ex plus -1ex minus
     -.2ex}{2.3ex plus .2ex}{\normalsize\bf}}
\def\subsection{\@startsection{subsection}{2}{\z@}{-3.25ex plus -1ex minus
 -.2ex}{1.5ex plus .2ex}{\normalsize\bf}}

\def\thebibliography#1{\section*{References\markboth
  {REFERENCES}{REFERENCES}}\list
  {[\arabic{enumi}]}{\settowidth\labelwidth{[#1]}\leftmargin\labelwidth
  \advance\leftmargin\labelsep
  \usecounter{enumi}}
  \def\newblock{\hskip .11em plus .33em minus -.07em}
  \sloppy
  \sfcode`\.=1000\relax}
 
\catcode `\@=12
\begin{document}

\vspace*{2.5cm}
\noindent
{ \bf Particles and Quantum Symmetries }\vspace{1.3cm}\\
\noindent
\hspace*{1in}
\begin{minipage}{13cm}

W{\l}adys{\l}aw Marcinek$^{1}$   \vspace{0.3cm}\\
 $^{1}$ Institute of Theoretical Physics, University
of Wroc{\l}aw\\
      \makebox[3mm]{ }Poland \\
\end{minipage}

\vspace*{0.5cm}

\begin{abstract}
\noindent
A system of interacting particles equipped with quantum
symmetry is described in an abstract algebraic way. The 
concept of quantum commutativity is used for description
of the algebra of quantum states of the system. The graded 
commutativity and particles in singular magnetic field is 
considered as an example. The generalized Pauli exclusion
principle is also mentioned.
\end{abstract}
\section{\hspace{-4mm}.\hspace{2mm} INTRODUCTION}
In this report particle systems intereacting with certain
external field and equipped with a quantum symmetry are 
considered. Our study is based on the assumption that the
system is characterized by a given Hopf algebra $H$ which 
act or coact on an unital and associative algebra $\ca$. 
The Hopf algebra $H$ describes quanta of external field 
and the algebra $\ca$ describes physical states of the 
system of particles. The action of $H$ on $\ca$ represents
the process of emission of quanta of the external field.
The coaction correspond for the process of absorption.
If the algebra $\ca$ is quantum commutative with respect 
to the given action or coaction of Hopf algebra $H$, 
then we say that the system posses a {\it quantum symmetry}. 
Note that quanum commutative algebras have been studied
previously by several authors, see \cite{sma,bae,mon,cowe} 
for example. In this paper we assume that the Hopf algebra 
$H$ is a group algebra $kG$, where $G$ is an abelian group 
and $k \equiv \com$ is the field of complex numbers. In
this particular case the coaction of $H$ on $\ca$ is
equivalent to the $G$--gradation of $\ca$ and the quantum 
commutativity is equivalent to the well--known graded 
commutativity \cite{mon,cowe}. 

A purely algebraic model of 
particles in singular magnetic field corresponding to 
the filling factor $v = \frac{1}{N}$ has been given by 
the author in \cite{top}. We describe here this model 
as an example of system with quantum symmetry. This 
symmetry is described as the gradation of the 
algebra $\ca$ by the group
$G \equiv Z_2 \oplus \cdots \oplus Z_2$ (N sumands).
We generalize our considerations for arbitrary fraction
$v = \frac{n}{N}$, where $n$ is the number of charged
particles per $N$ fluxes. 
Note that there is an approach in which the algebra 
$\ca$ play the role of noncommutative Fock space. 
In this attempt the creation and anihilation operators 
act on the algebra $\ca$ such that the creation operators 
act as the multiplication in $\ca$ and the anihilation 
ones act as a noncommutative contraction (noncommutative
partial derivatives) \cite{WM8,m10,ral}. One can use
such formalism in order to give the algebraic Fock
space representation \cite{wmq,com,sin}.
\section{\hspace{-4mm}.\hspace{2mm} MATHEMATICAL 
PRELIMINARIES}
Let us start with a brief review on the 
concept of quantum commutativity \cite{sma,cowe}.
It is well known that for a quasitriangular 
Hopf algebra $H$ there is the so--called $R$--matrix 
$R = \Sigma R_1 \ot R_2 \in H \ot H$ stisfing the 
quantum Yang--Baxter equation 
\be
(R \ot id)(id \ot R)(R \ot id) =
(id \ot R)(R \ot id)(id \ot R).
\ee
An algebra $\ca$ equipped with an action of 
$\tr : H \ot \ca \ra \ca$ is said to be 
{\it quantum commutative} if the following condition
\be
ab = \Sigma (R_2 \tr b)(R_1 \tr a)
\ee
is satisfied for every $a, b \in \ca$, see \cite{sma,bae}
for example. It is interesting that there is a dual notion 
of quantum commutativity. It is based on the coquasitriangular
Hopf algebras (CQHA) \cite{mon,cowe}. Let $H$ be a Hopf 
algebra over o field $k$. We use the following 
notation for the coproduct in $H$: if $h \in H$, 
then $\tg (h) := \Sigma h_1 \ot h_2 \in H \ot H$. 
A {\it coquasitriangular Hopf algebra} (a CQTHA) 
is a Hopf algebra $H$ equipped with a bilinear form 
$\lb -,-\rb : H \ot H \ra k$ such that
\be
\ba{l}
\Sigma \lb h_1, k_1 \rb k_2 h_2 = 
\Sigma h_1 k_1 \lb h_2, k_2 \rb ,\\
\lb h, kl \rb = \Sigma \lb h_1, k \rb \lb h_2, l \rb,\\
\lb hk, l\rb = \Sigma \lb h, l_2 \rb \lb k, l_1 \rb
\ea
\ee
for every $h, k, l \in H$. If such bilinear form $b$ exists
for a given Hopf algebra $H$, then we say that there is a 
{\it coquasitriangular structure} on $H$.

Let $H$ be a CQTHA with coquasitriangular structure given 
by a biform $\lb -,-\rb : H\ot H \ra k$ and let $\ca$ be 
a (right) $H$-comodule algebra with coaction $\rho$. 
Then the algebra $\ca$ is said to be {\it quantum commutative} 
with respect to the coaction of $(H, b)$ if an only if we
have the relation
\be
\ba{c}
a \ b = \Sigma \ \lb a_1, b_1 \rb \ b_0 \ a_0 ,
\ea
\ee
where $\rho (a) = \Sigma a_0 \ot a_1 \in \ca \ot H$, and
$\rho (b) = \Sigma b_0 \ot b_1 \in \ca \ot H$
for every $a, b \in \ca$, see \cite{cowe}. 
The Hopf algebra $H$ is said to be
a {\it quantum symmetry} for $\ca$.
Let $G$ be an arbitrary group, then the group algebra $H := kG$
is a Hopf algebra for which the comultiplication, the counit,
and the antypode are given by the formulae
$$
\ba{cccc}
\tg (g) := g \ot g,&\eta(g) := 1,&S(g) 
:= g^{-1}&\mbox{for} \ g \in G.
\ea
$$
respectively. If $H \equiv kG$, where $G$ is an abelian 
group, $k \equiv \com$ is the field of complex numbers, 
then the coquasitriangular structure on $H$ is given
as a bicharacter on $G$ \cite{mon}. Note that for
abelian groups we use the additive notation. A mapping 
$\ep : G \times G \ra \com \setminus \{0\}$ is said to 
be a {\it bicharacter} on $G$ if and only if we have 
the following relations
\be
\ep (\a, \b + \g) = \ep (\a, \b) \ep (\a , \g), \quad
\ep (\a + \b , \g) = \ep (\a , \g) \ep (\b , \g)
\ee
for $\a, \b, \g \in G$. If in addition
\be
\ep (\a, \b) \ep (\b, \a) = 1, 
\ee
for $\a, \b \in G$, then $\ep$ is said to be 
a {\it normalized bicharacter}. Note that this 
mapping are also said to be a commutation 
factor on $G$ and it has been studied previously by 
Scheunert \cite{sch}) in the context of color Lie algebras.
We restrict here our attention to normalized bicharacters 
only. It is interesting that the coaction of $H$ on certain 
space $E$, where $H \equiv kG$ is equivalent to the 
$G$-gradation of $E$, see \cite{mon}. In this case the 
quantum commutative algebra $\ca$ becomes graded commutative
\be
ab = \ep (\a , \b) b a
\ee
for homogeneous elements $a$ and $b$ of grade $\a$ and
$\b$, respectively.

Let $H$ be a CQTHA with coquasitriangular structure 
$\lb -,-\rb$. The family of all $H$-comodules forms 
a category $\cc = \cm^H$. The category $\cc$ is braided 
monoidal. The braid symmetry 
$\Psi \equiv \{\ps U, V, : U\ot V\ra V\ot U; U, V\in Ob\cc\}$
in $\cc$ is defined by the equation
\be
\ba{c}
\ps U, V, (u \ot v) = \Sigma\lb v_1 , u_1 \rb \ v_0 \ot u_0 ,
\label{coin}
\ea
\ee
where $\rho (u) = \Sigma u_0 \ot u_1 \in U \ot H$, and
$\rho (v) = \Sigma v_0 \ot v_1 \in V \ot H$ for every 
$u \in U , v \in V$. It is known that in an arbtrary 
braided monoidal category $\cc$ there is an algebra $\ca$ 
with braid commutative multiplication $m \circ \Psi = m$.
Obviously this algebra is quantum commutative.
\section{\hspace{-4mm}.\hspace{2mm} FUNDAMENTAL
ASSUMTIONS}
Let $E$ be a finite dimensional Hilbert space equipped with 
a basis $\th_i, i = 1,...,N = dim H$. We assume that $E$ is
a space of single particle quantum states for certain system of
particles. We also assume that there is a coaction $\rho_E$ 
of the Hopf algebra $H$ on the space $E$, i.e. a linear mapping 
$\rho_E : E \ra E \ot H$, which define a (right-) $H$-comodule 
structure on $E$. Then there is a braided monoidal category
$\cc \equiv \cc (E, \rho_H)$ which contains the field $\com$,
the space $E$, all tensors product and direct sums of these 
spaces, and some quotients. An algebra defined as the quotient 
$\ca \equiv \ca (E) := TE/I_{b}$, where $I_b$ is an 
ideal in $TE$ generated by elements of the form
\be
\ba{c}
u \ot v - \Sigma\lb v_1 , u_1 \rb \ v_0 \ot u_0
\label{alac}
\ea
\ee
and is said to be $\lb -,-\rb$-symmetric algebra over $E$, 
in the category $\cc$ \cite{mon}. The algebra $\ca$ 
describe partitions of particles dressed with quanta
of certain physical external field.

Now let us assume that $H = \com G$.
In this particular case the category $\cc = \cm^H$ 
becomes symmetric. The category $\cc = \cm^H$,
where $H := \com G$ for certain abelian group $G$ and 
$\lb -,-\rb$ is a bicharacter like above is denoted by 
$\cc \equiv \cc(G, \lb -,-\rb)$. Observe that the symmetry 
$\Psi$ in the formula (\ref{coin}) can be understood as 
a cointertwiner for corepresentations $\rho_U$ and $\rho_V$.
Note that if $E$ is a $H$-comodule, where $H = \com G$,
then $E$ is also a $G$-graded vector space, i.e 
$E = \bigoplus\limits_{\a \in G} \ E_{\a}$. 
This means that the above coaction is equivalent to 
$G$-gradation. The normalized bicharacter on the 
grading group $G$ has the following form
\be
\ba{c}
\ep (\a, \b ) = (-1)^{(\a |\b )} \ q^{<\a |\b >},
\label{don}
\ea
\ee
where $(-|-)$ is an integer-valued symmetric bi-form on $G$,
and $<-|->$ is a skew-symmetric integer-valued  bi-form on 
$G$, $q$ is some complex parameter \cite{zoz}. We use here 
the so-called standard gradation \cite{WM4}. This means that 
the grading group is $G \equiv Z^N := Z \oplus...\oplus Z$ 
($N$-sumands) for arbitrary $q$. If $q = exp(\frac{2\pi i}{n})$, 
$n \underline{>} 3$, then the grading group $G \equiv Z^N$ 
can be reduced to $G = Z_n \oplus...\oplus Z_n$. 
If $q = \pm 1$, then the grading group $G$ can be 
reduced to the group $Z_2 \oplus...\oplus Z_2$. 
In the standard gradation (see \cite{WM4}) the algebra
$\ca$ is generated by the relation
\be
\ba{c}
\th_{i} \ \th_{j} = \ep_{ij} \ \th_{j} \ \th_{i},
\label{laco}
\ea
\ee
where $\ep_{ij} := \ep(\si_i, \si_j)$ for $i,j = 1,...,N$, 
$\si_i$ is the $i$-th generator of the group $G$. In this
gradation the $i$--th generator $\th_i^a$ is a homogeneous 
element of grade $\si_i$ and describes a particle dressed 
with a single quantum. Here the quantum
commutativity with respect to the algebra $H \equiv \com G$ 
is the abovew graded commutativity \cite{sch,zoz,WM4}. 
Every element $\th_{\a}$ of the algebra $\ca$ can be given 
in the form
\be
\ba{c}
\th_{\a} = \th_1^{\a_1},..., \th_N^{\a_N}
\label{sta}
\ea
\ee
for $\a = \Sigma_{i=1}^N \ a_i \si_i$, $x_0 \equiv {\bf 1}$, 
where ${\bf 1}$ is the unit in $\ca$. 
\section{\hspace{-4mm}.\hspace{2mm}PARTICLES
IN SINGULAR MAGNETIC FIELD}
Now let us consider an algebraic model for a system of
charged particles moving in two-dimensional space $\cm$ 
with perpendicular magnetic field. Our fundamental 
assumption here is that the magnetic field is completely 
concentrated in vertical lines (fluxes). Every charged 
particle such as electron moving under influence of 
such singular magnetic field is transformed into a 
composite system which consists a charge and certain 
number of magnetic fluxes. We assume that there 
is in average $n$ electrons per $N$ fluxes. This 
means that the filling factor is $v = \frac{n}{N}$. 
We also assume that every magnetic flux can be
coupled to a charged particle such that the magnetic 
field of the flux is canceled. In this case we say that 
the flux is {\it bound} to a particle or {\it absorbed}
by it. The flux which is bound to a particle is said 
to be a {\it quasiparticle}. The particle is also said 
to be dressed by the flux. The flux not bound 
to a particle is said to be a {\it quasihole}.
Our system is characterized by all possible equivalence
classes of partitions of quasiparticle states and
quaisholes. We are going to construct a quantum
commutative algebra $\ca$ for description of such
partitions. First let us observe that the "effective" 
configuration space for the single particle is 
$\cm \equiv \real ^2 \setminus \{s_1,...,s_N\}$, 
where $s_1,...,s_N \in \real^2$ are points of 
intersection of magnetic lines with the plane. The 
fundamental group $\pi_1 (\cm)$ of $\cm$ is denoted 
by $G_N$. Let us denote by $\tilde{\si}_i$ the homotopy 
class of paths which arounds the singularity $s_i$ with 
the winding number equal to $1$ and not arounds any
other point $s_j , j \neq i$. The fundamental group 
$G_N$ is a free group generated by $\tilde{\si}_i$, 
$(i = 1,...,N)$ \cite{tqst}. We introduce an eqyuivalence 
relation in the loop space $\pi(\cm)$ as follows: If the phase
change of particle moving along a loop $\xi$ in $\cm$ is equal 
to the phase change of particle moving along an another loop 
$\xi'$ in $\cm$, then we say that such path are equivalent.
The fact that two lopps $\xi$ and $\xi'$ are equivalent 
is denoted by the expression $\xi \equiv_{ph} \xi'$. 
It is obvious that 
$\tilde{\si}_i\tilde{\si}_j\equiv_{ph}\tilde{\si}_j\tilde{\si}_i$
for $i \neq j$. This means that we have the relation
\be
(G_N/\equiv_{ph}) = Z \oplus...\oplus Z 
\quad \mbox{(N-sumands)}.
\ee
It is known that the coresponding group algebra $H := \com G$
is a coquasitriangular Hopf algebra. The coquasitriangular 
structure on $H$ is given by a bicharacter 
$\ep : \Gamma \times \Gamma \ra \com - \{0\}$ on $G$.
The normalized bicharacter $\ep$ for our model is given by 
the formula
\be
\ba{c}
\ep_{ij} = - (-1)^{N} (-1)^{\Omega_{ij}},
\label{com}
\ea
\ee
where $\Omega_{ij} = -\Omega_{ji} = 1$, for $i \neq j$, 
$\Omega_{ii} = 0$, $i, j = 1, 2,...,N$. This means that 
the grading group is
\be
\ba{c}
G \equiv Z_2^N := Z_2 \oplus ... \oplus Z_2\quad (N sumands)
\label{gru}
\ea
\ee
Now we describe the the algebra $\ca$ corresponding to our
model. In this case the algebra $\ca$ is generated by 
$\th_i^a (i=1,\ldots ,N; a=1,\ldots n)$ and relations
\be
\ba{lll}
\th_i^a \ \th_j^a =  \ep_{ij} \ \th_j^a \ \th_i^a&
\mbox{for} &a=b\\
\th_i^a \ \th_j^b = - \ep_{ij} \ \th_j^b \ \th_i^a&
\mbox{for} &a\neq b\\
(\th_i^a)^2 = 0&& ,
\ea
\ee
where $\ep$ is given by the formula (\ref{com}). 
The $G$--gradation of $\ca$ is given by the relation
\be
\mbox{grade} \th_i^a = \si_i ,
\ee
where $\si_i = (0\ldots1 \ldots 0)$ ($1$ on the $i$--th
place), is the $i$--th generator of the grading group 
$G = Z_2^N$. In our physical interpretation the generator 
$\th^a_i$ desribes the partition containing one 
particle dressing with the $i$--th flux (a quasiparticle)
and $N-1$ quasiholes. The statistics of 
quasiparticles is desrcibed by normalized bicharacter.
We also have here the generalized Pauli exclusion principle.
According to this principle partitions containing two or more
anticommuting quasiparticles must be excluded.
If $N$ is even, then we obtain the following relations 
\be
\ba{lll}
\th_i^a \ \th_j^a = - (-1)^{\Omega_{ij}} \ \th_j^a \ \th_i^a&
\mbox{for} &a=b\\
\th_i^a \ \th_j^b = (-1)^{\Omega_{ij}} \ \th_j^b \ \th_i^a&
\mbox{for} &a\neq b\\
(\th_i^a)^2 = 0&& .
\ea
\ee
For odd $N$ we obtain
\be
\ba{lll}
\th_i^a \ \th_j^a = (-1)^{\Omega_{ij}} \ \th_j^a \ \th_i^a&
\mbox{for} &a=b\\
\th_i^a \ \th_j^b = - (-1)^{\Omega_{ij}} \ \th_j^b \ \th_i^a&
\mbox{for} &a\neq b\\
(\th_i^a)^2 = 0&& 
\label{ant}
\ea
\ee
Let us consider a few examples.\\
{\bf Example 1.} Let us assume that $N = 2$ and $n = 1$. 
In this case the algebra $\ca$ is generating by two 
generators $\th_1$ and $\th_2$ satisfing the following 
relations
\be
\ba{l}
\th_1 \th_2 = \th_2 \th_1,\\
\th_1^2 = \th_2^2 = 0
\label{sqa}
\ea
\ee
and describes the systems containing in average one 
electron per two magnetic fluxes, i.e. the filling factor 
is $v = \frac{1}{2}$. Observe that there are two partitions
with one quasiparticle and one quasihole, namely $\th_1$
and $\th_2$. It is interesting that there is one partition 
which does not contain quasiholes. This is represent by 
the symmetric monomial $\th_1 \th_2$.\\
{\bf Example 2.} If $N = 3$ and $n = 1$, then the filling 
factor is $v = \frac{1}{3}$ and the algebra $\ca$ is 
genmerated by three generators $\th_1 , \th_2 , \th_3$ 
and relations
\be
\ba{lll}
\th_i \th_j = - \th_j \th_i,&\mbox{for}&i\neq j,\\
\th_1^2 = \th_2^2 = \th_3^2 = 0 .
\label{rsqa}
\ea
\ee
We have here three partitions $\th_1$, $\th_2$, and $\th_3$ 
which contain one quasiparticle and two quasiholes. Observe 
that the following partitions $\th_1 \th_2$, $\th_1 \th_3$,
$ \th_2 \th_3$ which contain two quasiparticles and one 
quasihole and the partition $\th_1 \th_2 \th_3$ which 
contains three quasiparticles are not physically 
admissible by the generalized Pauli exclusion principle. 
Hence in this case the single quasiparticle partitions
with two quasiholes can be physically realized!\\
{\bf Example 3.} For $N = 3$ and $n = 2$ we have the 
filling factor $v = \frac{2}{3}$ and the algebra $\ca$ is
genmerated by $\th_i^a$, where $a=1, 2$ and $i=1, 2, 3$
such that
\be
\ba{lll}
\th_i^a \th_j^a = - \th_j^a \th_i^a,&\mbox{for}&a=b,\\
\th_i^a \th_j^b = \th_j^b \th_i^a,&\mbox{for}&a\neq b,
\label{sqap}
\ea
\ee
for $i \neq j$, $(\th_i^a)^2 = 0$ and
\be
\ba{lll}
\th_i^a \th_j^a = \th_j^a \th_i^a,&\mbox{for}&a=b,\\
\th_i^a \th_j^b = - \th_j^b \th_i^a,&\mbox{for}&a\neq b,
\label{tqap}
\ea
\ee
for $i = j$. Observe that quasiparticles $\th_i^a$ and 
$\th_j^b$ for $i \neq j$ and $a \neq b$ commutes and 
partitions with two quasiholes are admissible.

\vspace{\baselineskip}

\end{document}